\documentstyle[11pt]{article}
\title{Examples illustrating some aspects of the weak Deligne-Simpson problem
\footnote{Research partially supported 
by INTAS grant 97-1644}}
\author{Vladimir Petrov Kostov\\ \\ \hspace{7cm}
{\sl To the memory of my mother}} 
\date{}
\newtheorem{tm}{Theorem}
\newtheorem{lm}[tm]{Lemma}

\newtheorem{prop}[tm]{Proposition}
\newtheorem{rem}[tm]{Remark}
\newtheorem{rems}[tm]{Remarks}
\newtheorem{defi}[tm]{Definition}
\newtheorem{ex}[tm]{Example}

\newtheorem{nota}[tm]{Notation}
\begin{document}
\maketitle 

\begin{abstract}
We consider the variety of $(p+1)$-tuples of matrices $A_j$ (resp. $M_j$) 
from given conjugacy classes $c_j\subset gl(n,{\bf C})$ 
(resp. $C_j\subset GL(n,{\bf C})$) such that $A_1+\ldots +A_{p+1}=0$ 
(resp. $M_1\ldots M_{p+1}=I$). This variety is connected with the 
weak {\em Deligne-Simpson problem: give necessary and sufficient conditions on 
the choice of the conjugacy classes $c_j\subset gl(n,{\bf C})$ 
(resp. $C_j\subset GL(n,{\bf C})$) so that there 
exist $(p+1)$-tuples with trivial centralizers of matrices $A_j\in c_j$ 
(resp. $M_j\in C_j$) whose sum equals 0 (resp. whose product equals $I$).} 
The matrices $A_j$ (resp. $M_j$) are interpreted as matrices-residua of 
Fuchsian linear systems (resp. as 
monodromy operators of regular linear systems) on Riemann's sphere. We 
consider examples of such varieties of   
dimension higher than the expected one due to the presence 
of $(p+1)$-tuples with non-trivial centralizers; in one of the examples 
the difference between the two dimensions is $O(n)$.\\  


{\bf AMS classification index:} 15A30, 20G05
\end{abstract}

\section{Introduction}
\subsection{Formulation of the (weak) Deligne-Simpson problem
\protect\label{WDSP}}

In the present article we consider examples   related to the 
{\em Deligne-Simpson problem (DSP)}. The problem stems from the analytic 
theory of linear systems of ordinary differential equations but its 
formulation is purely algebraic: 

{\em Give necessary and sufficient conditions on the choice of the $p+1$ 
conjugacy classes $c_j\subset gl(n,{\bf C})$, resp. 
$C_j\subset GL(n,{\bf C})$, so that there exist irreducible $(p+1)$-tuples of 
matrices $A_j\in c_j$, $A_1+\ldots +A_{p+1}=0$, resp. of matrices 
$M_j\in C_j$, $M_1\ldots M_{p+1}=I$.}

Here $I$ stands for the identity matrix and ``irreducible'' 
means ``with no non-trivial common invariant subspace''. 
The version with matrices $A_j$ (resp. $M_j$) is called the {\em additive} 
(resp. the {\em multiplicative}) one. 
The matrices $A_j$ are interpreted as {\em matrices-residua} of 
{\em Fuchsian} systems on Riemann's sphere (i.e. linear systems of ordinary 
differential equations with logarithmic poles). The sum of all 
matrices-residua of a Fuchsian system equals 0.  

The matrices $M_j$ 
are interpreted as 
{\em monodromy operators} 
of meromorphic linear  
{\em regular} systems on Riemann's sphere (i.e. linear systems of ordinary 
differential equations with moderate growth rate of the solutions at the 
poles). (Fuchsian systems are always regular.) A {\em monodromy operator} 
of a regular system is a linear operator acting on its solution space which 
maps the solution with a given initial value at a given {\em base point} $a_0$ 
onto the value at $a_0$ of its analytic continuation along some closed 
contour. 

The monodromy operators generate the 
{\em monodromy group}. One usually chooses as generators of the monodromy 
group operators defined by contours which are freely homotopic to small 
loops each circumventing counterclockwise one of the poles of the system. 
For a suitable indexation of the poles the product of these generators equals 
$I$ (and this is the only relation which they a priori satisfy). 

\begin{rem}
In the multiplicative version  
the classes $C_j$ are interpreted as 
{\em local monodromies} around the poles and the DSP admits the 
following interpretation: 

For what $(p+1)$-tuples of local monodromies do there exist irreducible  
monodromy groups with such local monodromies.
\end{rem}

The monodromy group of a regular system is the only 
invariant of a regular system under the linear changes of the dependent 
variables meromorphically depending on the time. Therefore the multiplicative 
version is more important than the additive one; nevertheless, the additive 
one is easier to deal with when computations are to be performed and one can 
easily deduce corollaries concerning the multiplicative version as well 
due to Remark~\ref{resmon}. 

\begin{rem}\label{resmon}
If $A$ denotes a matrix-residuum at a given pole of a Fuchsian system 
and if $M$ denotes the 
corresponding operator of local monodromy, then 
in the absence of non-zero integer differences 
between the eigenvalues of $A$ the operator $M$ is conjugate to 
$\exp (2\pi iA)$. 
\end{rem}

By definition, the {\em weak DSP} is the 
DSP in which instead of 
irreducibility of the 
$(p+1)$-tuple of matrices one requires only its centralizer to be trivial. 
We say that the DSP (resp. the weak DSP) is {\em solvable} for a given 
$(p+1)$-tuple of 
conjugacy classes $c_j$ or $C_j$ if there exist matrices $A_j\in c_j$ whose 
sum is 0 or matrices $M_j\in C_j$ whose product is $I$ such that their 
$(p+1)$-tuple is irreducible (resp. with trivial centralizer). By definition, 
the (weak) DSP is solvable for $n=1$.

For given conjugacy classes $c_j$  
satisfying the condition $\sum$Tr$(c_j)=0$ 
and also the ones of Theorem~\ref{necessary} below consider the variety 

\[ {\cal V}(c_1,\ldots ,c_{p+1})=\{ (A_1,\ldots ,A_{p+1})~|~A_j\in c_j~,~
A_1+\ldots +A_{p+1}=0\} ~. \]
(One can define such a variety in a similar way in the case of matrices $M_j$ 
as well.) 

If the eigenvalues are generic (see the precise definition in 
Subsection~\ref{Gen}), then the variety 
${\cal V}(c_1,\ldots ,c_{p+1})$ (or just ${\cal V}$ for short) is smooth, 
see \cite{Ko3}. If not, then it can have a complicated stratified 
structure defined by the 
invariants of the $(p+1)$-tuple of matrices; the 
centralizers of the 
$(p+1)$-tuples might be trivial on some strata and non-trivial on others; 
finally, a stratum on which the 
centralizer is non-trivial can be of greater dimension than the one of a 
stratum on which it is trivial; see \cite{Ko3} for some examples. 

The aim of the present paper is to give further examples of varieties 
${\cal V}$ and to discuss their stratified structure and dimension of the 
strata. This will be explained in some more detail in Subsection~\ref{rigind} 
after some necessary notions will be introduced in the next two subsections.

\begin{rem}
In what follows the sum of the matrices $A_j$ is always presumed to be 
0 and the product of the matrices $M_j$ is always presumed to be $I$.
\end{rem}

\begin{nota}
Double subscripts indicate matrix entries. We denote by $E_{i,j}$ the matrix 
having zeros everywhere except in position $(i,j)$ where it has a unit.
\end{nota} 

\subsection{Necessary conditions for the solvability of the (weak) DSP}

The known results concerning the (weak) DSP are exposed in \cite{Ko1}. We 
recall in this and in the next subsection only the most necessary ones.

\begin{defi}
A {\em Jordan normal form (JNF) of size $n$} is a collection of positive 
integers 
$\{ b_{i,l}\}$ whose sum is $n$ 
where $b_{i,l}$ is the size of the $i$-th Jordan block with 
the $l$-th eigenvalue; the eigenvalues are presumed distinct and for $l$ fixed 
the numbers $b_{i,l}$ form a non-increasing sequence. Denote by $J(C)$ (resp. 
$J(A)$) the 
JNF defined by the conjugacy class $C$ (resp. by the matrix $A$).  
\end{defi}

\begin{defi}
For a conjugacy class $C$ in $GL(n,{\bf C})$ or $gl(n,{\bf C})$ denote by 
$d(C)$ its dimension; recall that it is always even. For a matrix $Y\in C$ set 
$r(C):=\min _{\lambda \in {\bf C}}{\rm rank}(Y-\lambda I)$. The integer 
$n-r(C)$ is the maximal number of Jordan blocks of $J(Y)$ with one and the 
same eigenvalue. Set $d_j:=d(C_j)$ (resp. $d(c_j)$), $r_j:=r(C_j)$ 
(resp. $r(c_j)$). The quantities 
$r(C)$ and $d(C)$ depend only on the JNF $J(Y)=J^n$, not 
on the eigenvalues, so we write sometimes $r(J^n)$ and $d(J^n)$. 
\end{defi}

The following two inequalities are necessary conditions for the existence 
of irreducible $(p+1)$-tuples of matrices $A_j$ or $M_j$ (their 
necessity in the multiplicative version was proved by C.Simpson, 
see \cite{Si}, and in the additive one by the author, see \cite{Ko2}):

\[ d_1+\ldots +d_{p+1}\geq 2n^2-2~~~~~(\alpha _n)~~,~~~~~
{\rm for~all~}j,~r_1+\ldots +\hat{r}_j+\ldots +r_{p+1}\geq n~~~~~
(\beta _n)~~~.\]

The inequality 

\[ r_1+\ldots +r_{p+1}\geq 2n~~~~~(\omega _n)\]
is not a necessary condition (note that it implies $(\beta _n)$) but it is 
``almost sufficient'', i.e. sufficient in most part of the cases, 
see the details in \cite{Ko1}. 
 
We formulate below a necessary condition for the solvability of the (weak) 
DSP which is 
a condition upon the $p+1$ JNFs $J_j^n=J(c_j)$ or $J(C_j)$ 
($j=1,\ldots ,p+1$, the upper index 
indicates the size of the matrices) but not upon the classes $c_j$ or $C_j$ 
themselves. 

\begin{defi}
For a given $(p+1)$-tuple of JNFs $J_j^n$ with $n>1$, which satisfies 
condition $(\beta _n)$ and doesn't satisfy condition 
$(\omega _n)$ set $n_1=r_1+\ldots +r_{p+1}-n$. Hence, $n_1<n$ and 
$n-n_1\leq n-r_j$. Define 
the $(p+1)$-tuple of JNFs $J_j^{n_1}$ as follows: to obtain the JNF 
$J_j^{n_1}$ 
from $J_j^n$ one chooses one of the eigenvalues of $J_j^n$ with 
greatest number $n-r_j$ of Jordan blocks, then decreases  
by 1 the sizes of the $n-n_1$ {\em smallest} Jordan blocks with this 
eigenvalue and deletes the Jordan blocks of size 0. Denote this construction 
by $\Psi :(J_1^n,\ldots ,J_{p+1}^n)\mapsto (J_1^{n_1},\ldots ,J_{p+1}^{n_1})$ 
or just by $\Psi$ for short.
\end{defi} 

\begin{tm}\label{necessary}
If the (weak) DSP is solvable for a given $(p+1)$-tuple of conjugacy classes 
$c_j$ or $C_j$ defining the JNFs $J_j^n$, satisfying condition 
$(\beta _n)$ and not satisfying condition $(\omega _n)$, then the 
map $\Psi$ iterated as long as defined stops  
at a $(p+1)$-tuple of JNFs $J_j^{n'}$ either 
satisfying condition $(\omega _{n'})$ or with $n'=1$.
\end{tm} 

The theorem can be deduced from \cite{Ko1}, see Theorem 8 there. 

\begin{rem}
One can show that the results formulated by means of the map $\Psi$ do not 
depend on the choice of an eigenvalue with maximal number of Jordan blocks 
belonging to it whenever such a choice is possible.
\end{rem}

\subsection{Generic eigenvalues and (poly)multiplicity vectors
\protect\label{Gen}}

We presume in the case of matrices $M_j$ 
the necessary condition $\prod \det (C_j)=1$ to hold. In the case of 
matrices $A_j$ this is the condition 
$\sum$Tr$(c_j)=0$. In terms of the eigenvalues $\sigma _{k,j}$   
(resp. $\lambda _{k,j}$) of the matrices from $C_j$ (resp. $c_j$) repeated 
with their multiplicities, this condition reads     
$\prod _{k=1}^n\prod _{j=1}^{p+1}\sigma _{k,j}=1$  
(resp. $\sum _{k=1}^n\sum _{j=1}^{p+1}\lambda _{k,j}=0$). 

\begin{defi}
An equality of the form  
$\prod _{j=1}^{p+1}\prod _{k\in \Phi _j}\sigma _{k,j}=1$, resp. 
$\sum _{j=1}^{p+1}\sum _{k\in \Phi _j}\lambda _{k,j}=0$, is called a 
{\em non-genericity relation};  
the sets $\Phi _j$ contain one and the same number $<n$ of indices  
for all $j$. Eigenvalues satisfying none of these relations are called 
{\em generic}. If one replaces for all $j$ the sets $\Phi _j$ by their 
complements in $\{ 1,\ldots, n\}$, then one obtains another non-genericity 
relation which we identify with the initial one.
\end{defi}

\begin{rems}
1) Reducible $(p+1)$-tuples exist only for non-generic eigenvalues. Indeed, if 
the $(p+1)$-tuple is block upper-triangular, then the eigenvalues of each 
diagonal block satisfy some non-genericity relation. 

2) For generic eigenvalues the conditions of Theorem~\ref{necessary} are 
sufficient as well, see \cite{Ko1} (Theorem~8), \cite{Ko2} and \cite{Ko6}.
\end{rems}

\begin{rem}\label{nongeneric}
Condition $(\beta _n)$ admits the following generalizations which in certain  
cases of non-generic eigenvalues are stronger than $(\beta _n)$ itself -- 
these are the inequalities  

\[ \min _{b_j\in {\bf C}, b_1+\ldots +b_{p+1}=0}
\sum _{j=1}^{p+1}{\rm rk}(A_j-b_jI)\geq 2n~~,~~\min _{b_j\in {\bf C}^*, 
b_1\ldots b_{p+1}=1}
\sum _{j=1}^{p+1}{\rm rk}(b_jM_j-I)\geq 2n~~~~~~~(\delta _n).\]
which are necessary conditions for the existence of irreducible 
$(p+1)$-tuples 
of matrices $A_j$ or $M_j$ (see \cite{Ko1}, Lemma 10 and the line after it). 
\end{rem}

\begin{defi}
A {\em multiplicity vector (MV)} is a vector with positive integer components 
whose sum is $n$ and 
which are the multiplicities of the eigenvalues of an $n\times n$-matrix. 
In the case of diagonalizable matrices the MV defines 
completely the JNF. A {\em polymultiplicity vector (PMV)} is a $(p+1)$-tuple 
of multiplicity vectors, the ones of the eigenvalues of the matrices $A_j$ or 
$M_j$. 
\end{defi}

\begin{rem}\label{rd}
For a diagonal JNF $J^n$ defined by the MV $(m_1,\ldots ,m_s)$, 
$m_1\geq \ldots \geq m_s$, one has $r(J^n)=m_2+\ldots +m_s$ and 
$d(J^n)=n^2-m_1^2-\ldots -m_s^2$. If the  
JNF $J(c_j)$ or $J(C_j)$ is diagonal, then the construction $\Psi$ (see the 
previous subsection) results in decreasing the greatest component of the 
$j$-th MV by $n-n_1$.
\end{rem}

\subsection{The index of rigidity and the expected dimension of the 
variety ${\cal V}$\protect\label{rigind}}

\begin{defi}
Call {\em index of rigidity} of the $(p+1)$-tuple of  
conjugacy classes $c_j$ or $C_j$ (or of the $(p+1)$-tuple of JNFs 
defined by them) the quantity 
$\kappa =2n^2-d_1-\ldots -d_{p+1}$. This notion was introduced in 
\cite{Ka}.
\end{defi}

\begin{rems}\label{Gabber}
1) If condition $(\alpha _n)$ holds, then $\kappa$ can take the values 
$2$, $0$, $-2$, $-4$, $\dots$.

2) If $\kappa =2$ and the DSP is solvable for given conjugacy classes, 
then such $(p+1)$-tuples are unique up 
to conjugacy, see \cite{Ka} and \cite{Si} for the multiplicative version; from 
this result one easily deduces the uniqueness in the additive version.

3) For $\kappa =2$ the coexistence of irreducible and reducible $(p+1)$-tuples 
of matrices $M_j$ is impossible, see \cite{Ka}, Theorem~1.1.2, or \cite{Ko3}, 
Theorem~18. One can easily deduce from this fact that the same is true for 
matrices $A_j$.
\end{rems}

Recall that the variety ${\cal V}$ was defined in Subsection~\ref{WDSP}. 

\begin{rems}\label{expected}
1) If ${\cal V}$ is nonempty and if the eigenvalues are generic, then it 
contains only irreducible $(p+1)$-tuples, it is smooth and its 
dimension equals $1-\kappa +n^2$, see \cite{Ko3}. 

2) If on some stratum of ${\cal V}$ the 
centralizer is trivial, then the stratum is smooth and 
its dimension equals $1-\kappa +n^2$; we call 
this dimension the {\em expected dimension} of ${\cal V}$, see \cite{Ko3}, 
Proposition~2.
\end{rems}

In the present paper we consider examples of varieties ${\cal V}$ for 
conjugacy classes satisfying the conditions of Theorem~\ref{necessary}. The 
first of them (see Section~\ref{fex}) is 
with $\kappa =2$, $n$ odd and $p=3$. We discuss its stratified 
structure and we show 
that it contains a stratum (on which the centralizer is non-trivial) 
of dimension $n^2+(n-3)/2$, i.e. exceeding the 
expected one by $(n-1)/2$, as well as strata of dimensions 
$n^2+s-1$ for $s=1,2,\ldots ,(n-1)/2$. And it contains strata of  
expected dimension $n^2-1$ on which the centralizer is trivial. 
  
The second example (see Section~\ref{secex}) is one with $\kappa =2$ where the 
variety ${\cal V}$ is not connected (but its closure is).

In Section~\ref{tex} we show that for $(-\kappa )$ arbitrarily high there 
exist examples of varieties ${\cal V}$ in which certain strata where the 
centralizer is non-trivial do not belong to the closures of the strata where 
it is trivial and are of dimension higher than the expected one. This provides 
a negative answer to a question stated in \cite{Ko3}.

\section{An example with $p=3$\protect\label{fex}}
\subsection{Description of the example}

Consider for $p=3$, $n=2k+1$, $k\in {\bf N}^*$, the PMV 
$\Lambda (k)=((k+1,k),(k+1,k),(k+1,k),(k+1,k))$ (the matrices $A_j$ are 
presumed diagonalizable; a similar example can be given for matrices $M_j$ as 
well). Denote the respective 
eigenvalues of the matrices $A_j$ 
with such a PMV by $\lambda _j$, $\mu _j$ ($\lambda _j$ is of multiplicity 
$k+1$, one has $\lambda _j\neq \mu _j$, $j=1,2,3,4,$). 
One has $d_j=2k(k+1)$, 
see Remark~\ref{rd}, hence, $\kappa =2$. The PMV $\Lambda (k)$ satisfies 
the conditions of Theorem~\ref{necessary} -- one has 
$\Psi (\Lambda (k))=\Lambda (k-1)$ (see Remark~\ref{rd}) 
and the iterations of $\Psi$ stop at a quadruple of JNFs of size 1. 

We assume that 
$\sum _{j=1}^4\lambda _j=0~(A)$ is the 
only non-genericity relation satisfied by the eigenvalues (note that 
it implies $\sum _{j=1}^4\mu _j=0$). 

\begin{prop}\label{oneortwo}
Any quadruple of matrices $A_j$ like above whose sum is 0 is up to conjugacy 
block upper-triangular, with diagonal blocks of sizes 1 or 2. The diagonal 
blocks of size 1 equal either $(\lambda _1,\lambda _2,\lambda _3,\lambda _4)$ 
or $(\mu _1,\mu _2,\mu _3,\mu _4)$. The restriction of $A_j$ to a diagonal 
block of size 2 has eigenvalues $\lambda _j$, $\mu _j$.
\end{prop}

The propositions from this subsection are proved in the next ones.
   
\begin{ex}\label{exB}
There exist irreducible quadruples of $2\times 2$-matrices $B_j$ 
whose sum is 0 and with eigenvalues $\lambda _j$, $\mu _j$:

\[ B_1=\left( \begin{array}{cc}\lambda _1&1\\0&\mu _1\end{array}\right) ~,~
B_2=\left( \begin{array}{cc}\lambda _2&-1\\0&\mu _2\end{array}\right) ~,~
B_3=\left( \begin{array}{cc}\lambda _3&0\\u&\mu _3\end{array}\right) ~,~
B_4=\left( \begin{array}{cc}\lambda _4&0\\-u&\mu _4\end{array}\right) \]
where $u\in {\bf C}^*$. 
\end{ex}

\begin{prop}\label{Pi}
1) The variety $\Pi$ 
of quadruples of diagonalizable $2\times 2$-matrices $B_j$ with 
eigenvalues $\lambda _j$, $\mu _j$ and such that $B_1+\ldots +B_4=0$ is 
connected. 

2) Its subvariety $\Pi _0$ consisting of all such irreducible quadruples is 
also connected.
\end{prop}

\begin{ex}\label{exH}
For $l\in {\bf N}^*$ there exist upper-triangular quadruples of 
$(2l+1)\times (2l+1)$-matrices $H_j$ with zero sum, with trivial centralizers 
and with eigenvalues $\lambda _j$, $\mu _j$ of multiplicity $l+1$, $l$ 
(the matrix $I$ is $l\times l$):

\[ H_1=\left( \begin{array}{ccc}\lambda _1&0&0\\0&\lambda _1I&0\\0&0&\mu _1I
\end{array}\right) ~,~
H_2=\left( \begin{array}{ccc}\lambda _2&0&0\\0&\lambda _2I&I\\0&0&\mu _2I
\end{array}\right) ~,~
H_3=\left( \begin{array}{ccc}\lambda _3I&0&I\\0&\lambda _3&0\\0&0&\mu _3I
\end{array}\right) \]
(pay attention to the block-decomposition of $H_3$ which is different from 
the one of $H_1$ and $H_2$). If for $Z\in gl(2l+1,{\bf C})$ one has 
$[Z,H_1]=[Z,H_2]=[Z,H_3]=0$, then $[Z,H_1]=0$ implies that $Z$ is 
block-diagonal, with diagonal blocks of sizes $(l+1)\times (l+1)$ and 
$l\times l$. One deduces then from $[Z,H_2]=[Z,H_3]=0$ that $Z$ is scalar 
(the details are left for the reader).
\end{ex}

\begin{rem}\label{otherkinds}
Examples similar to the above one 
can be given for other permutations of the eigenvalues on the 
diagonal as well 
(e.g. when all $\mu _j$ come first followed by all $\lambda _j$). For some 
permutations there exist no examples of such upper-triangular quadruples 
with trivial centralizers (e.g. when the first and the last eigenvalues on 
the diagonal are equal -- in this case the matrices commute with $E_{1,2l+1}$).
\end{rem}

\begin{defi}
An irreducible quadruple of diagonalizable $2\times 2$-matrices with 
eigenvalues $(\lambda _j, \mu _j)$ whose sum is 0 is said to be of 
{\em type $B$}. (Example~\ref{exB} shows such a quadruple.) 

An upper-triangular quadruple with trivial centralizer of diagonalizable 
$h\times h$-matrices ($h=2l+1$) with eigenvalues $\lambda _j$, $\mu _j$ of 
multiplicities 
$l+1$, $l$ is said to be of {\em type $H_h$}. 
(Example~\ref{exH} shows such a quadruple.) 

A stratum of ${\cal V}$ the quadruples of which up to conjugacy are 
block-diagonal, with $s$ diagonal blocks of type $B$ defining  
non-equivalent representations and with one diagonal 
block of type $H_{n-2s}$ is said to be of {\em type $HB_s$}.
\end{defi}

\begin{prop}\label{HBs}
1) A stratum of 
type $HB_s$ is locally a smooth algebraic variety of dimension 
$n^2+s-1$. 

2) It is globally connected. 
\end{prop}

\begin{prop}\label{connected} 
The variety ${\cal V}$ from the example is connected.
\end{prop}

{\bf Conclusive remarks.} 
As we saw, the stratum of type $HB_0$ (Example~\ref{exH}) 
consists of quadruples with trivial 
centralizers and is of dimension $n^2-1$ (the expected one) 
while the one of type $HB_{(n-1)/2}$ 
is of dimension $n^2+(n-3)/2$. All intermediate dimensions are attained 
on the strata $HB_s$, see Proposition~\ref{HBs}. For $s>0$ they consist of 
quadruples with non-trivial centralizers (they are block-diagonal up to 
conjugacy). Except the strata of type $HB_s$ there are other 
strata of ${\cal V}$ with non-trivial centralizer, e.g. such on which the 
representation is a direct sum of some representations of type $B$ and a 
representation with a non-trivial centralizer as mentioned in 
Remark~\ref{otherkinds}.

\subsection{Proof of Proposition~\protect\ref{oneortwo}}

$1^0$. It is clear that there exist diagonal quadruples of matrices $A_j$ 
like above 
whose sum is 0 (their first $k+1$ diagonal entries equal $\lambda _j$ and the 
last $k$ ones equal $\mu _j$, see the non-genericity relation $(A)$). It 
follows from 3) of Remarks~\ref{Gabber} that 
there exist no irreducible such quadruples. 

$2^0$. Any reducible quadruple can be 
conjugated to a block upper-triangular form. The restriction of the quadruple 
to each diagonal block $B$ is presumed to define an irreducible 
representation. 

If the size $l$ of the block is odd and $>1$, then the minimal possible 
value of $\kappa$ for this block is $2$ and it is attained only when for 
each $j$ the 
multiplicities of $\lambda _j$ and $\mu _j$ as eigenvalues of $A_j|_B$ 
equal $(s+1,s)$ or $(s,s+1)$ where $l=2s+1$. (To prove this one can 
use Remark~\ref{rd}.) The absence of non-genericity 
relations other than $(A)$ implies that the multiplicity of $\lambda _j$ 
(and, hence, the one of $\mu _j$) is one and the same for all $j$. 
However, the existence of 
diagonal quadruples of matrices $A_j$ of size $l$ with such multiplicities 
of $\lambda _j, \mu _j$ implies that such blocks $B$ 
do not exist.

$3^0$. If the size of $B$ is $2m$, $m\in {\bf N}^*$, then the minimal possible 
value of $\kappa$ is 0 and it is attained only when the multiplicities of 
$\lambda _j$ and $\mu _j$ as eigenvalues of $A_j|_B$ equal $(m,m)$ (the easy 
computation is left for the reader). Such 
blocks $B$ exist only for $m=1$, see \cite{Ko5}. 

$4^0$. If the size of $B$ is $2m$, $m\in {\bf N}^*$, and if $\kappa =2$, then 
this can happen only if for three of the indices $j$ the multiplicities of 
$\lambda _j$ and $\mu _j$ as eigenvalues of $A_j|_B$ equal $(m,m)$ and for the 
fourth one (say, for $j=4$) 
they equal $(m-1,m+1)$ or $(m+1,m-1)$ (we leave the proof for the 
reader again). This together with $(A)$ implies that $\lambda _4=\mu _4$ 
which is impossible. Hence, such blocks $B$ do not exist. 

Hence, only blocks of size 1 and of size 2  
are possible to occur on the diagonal (for the ones of size 2 see $3^0$). The 
fact that $(A)$ is the only non-genericity relation implies that the blocks 
of size 1 equal either $(\lambda _1,\lambda _2,\lambda _3,\lambda _4)$ 
or $(\mu _1,\mu _2,\mu _3,\mu _4)$.

The proposition is proved.~~~~~$\Box$

\subsection{Proof of Proposition~\protect\ref{Pi}} 

$1^0$. Prove 1). Denote by $c_j^*$ the conjugacy class of the matrix $B_j$.  
Denote by $\tau$ the quantity tr$(B_1+B_2)$. By varying the matrices $B_1$ 
and $B_2$ (resp. $B_3$ and $B_4$) within their conjugacy classes 
one can obtain as their sum $S=B_1+B_2$ 
(resp. as $-(B_3+B_4)$) any non-scalar matrix from the 
set $\Delta (\tau )$ of $2\times 2$-matrices with trace 
equal to $\tau$. 

Indeed, if $S_{1,2}=g\neq 0$, then set 
$B_1=\left( \begin{array}{cc}\lambda _1&0\\u&\mu _1\end{array}\right)$, 
$B_2=\left( \begin{array}{cc}h&g\\w&\lambda _2+\mu _2-h\end{array}\right)$. 
One fixes first $h$ to obtain the necessary entry  
$S_{1,1}$. 
One has $g\neq 0$, hence, there exists a unique $w$ satisfying the 
condition $\det (B_2)=\lambda _2\mu _2$; after this one chooses $u$ to obtain 
the necessary entry $S_{2,1}$. 

$2^0$. If $S_{1,2}=0$, then one can conjugate $S$ by some matrix 
$Y\in GL(2,{\bf C})$ to obtain the condition 
$S_{1,2}\neq 0$, find the matrices $B_1$ and $B_2$ like above 
and then conjugate them (and $S$) by 
$Y^{-1}$. This is possible to do because $S$ is not scalar.

The sets $\Delta (\tau )$ and $\Delta (\tau )\backslash \{ \tau I/2\}$  
being connected so is the 
variety $\Pi$. Indeed, one has  
$\Pi =\{ (B_1,B_2,B_3,B_4)|B_j\in c_j^*,B_1+B_2=-(B_3+B_4)\}$.

$3^0$. Prove 2). If the quadruple of matrices $B_j$ is reducible, then so is 
the couple $B_1,B_2$, hence, the eigenvalues 
of $B_1+B_2$ equal either $(\lambda _1+\lambda _2,\mu _1+\mu _2)$ or 
$(\lambda _1+\mu _2,\mu _1+\lambda _2)$. 

The subset $\Delta ^0(\tau )$ 
of $\Delta (\tau )$ defined by the condition the eigenvalues of a matrix from 
$\Delta (\tau )$ to equal either $(\lambda _1+\lambda _2,\mu _1+\mu _2)$ or 
$(\lambda _1+\mu _2,\mu _1+\lambda _2)$ is a proper subvariety of the 
smooth irreducible variety $\Delta (\tau )$. Therefore the connectedness of 
$\Pi _0$ is proved just like the one of $\Pi$, by replacing $\Delta (\tau )$ 
by $\Delta (\tau )\backslash \Delta ^0(\tau )$.~~~~$\Box$

\subsection{Proof of Proposition~\protect\ref{HBs}}

$1^0$. Prove 1). 
Denote by $\Sigma$ the variety of block-diagonal quadruples whose first $s$ 
diagonal blocks are of type $B$ and the last block {\em up to conjugacy} is 
of type $H_{n-2s}$. Each of the 
first $s$ blocks defines a smooth 
variety of dimension 5, see Remarks~\ref{expected}. The last block 
defines a smooth variety of dimension $(n-2s)^2-1$. Hence, 
dim$\Sigma =5s+(n-2s)^2-1$. To deduce dim$(HB_s)$ from dim$\Sigma$ one has to 
add to dim$\Sigma$ the dimension of a transversal $T$ to the group of 
infinitesimal conjugations preserving the block-diagonal 
form of the quadruple. This is the group ${\cal G}$ 
of block-diagonal matrices (which are deformations of $I$) with 
the same sizes of the diagonal blocks as the ones of the quadruple (we leave 
the proof of this statement for the reader; use the fact that  
the diagonal blocks define non-equivalent representations the first $s$ 
of which of type $B$ and the last of type $H_{n-2s}$). Hence, 
dim${\cal G}=4s+(n-2s)^2$,  dim$T=n^2-4s-(n-2s)^2$ and 
dim${\cal V}=n^2-4s-(n-2s)^2+5s+(n-2s)^2-1=n^2+s-1$.

The stratum of type $HB_s$ is locally diffeomorphic to $\Sigma \times T$, 
hence, it is smooth.

$2^0$. Prove 2). The variety $\Pi _0$ of quadruples 
of matrices of type $B$ is connected, see Proposition~\ref{Pi}. It is smooth 
as well, hence, it is irreducible. 
Hence, the cartesian product of $s$ copies of $\Pi _0$ is connected; if one 
deletes from it the subvariety on which two of the representations are 
equivalent, then the resulting variety is still connected. 

Hence, the variety $\Sigma$ is connected. The connectedness of $\Sigma$ and 
the one of $GL(n,{\bf C})$ imply the one of the stratum of type 
$HB_s$.~~~~~$\Box$

\subsection{Proof of Proposition~\protect\ref{connected}}

Every quadruple of matrices $A_j$ from ${\cal V}$ can be conjugated 
to a block upper-triangular form with diagonal blocks of sizes 1 or 2 
(Proposition~\ref{oneortwo}). Conjugate the quadruple by a suitable 
one-parameter family 
of diagonal matrices to make all entries above the diagonal tend to 0 while 
preserving the diagonal blocks. The limit quadruple (denoted by 
$(A_1',\ldots ,A_4')$) also belongs to 
${\cal V}$. Indeed, the restriction of $A_j'$ to each diagonal block of size 2 
is diagonalizable, hence, $A_j'$ is diagonalizable, the eigenvalues and their 
multiplicities are the same as for $A_j$ and the sum of the matrices $A_j'$ 
is 0. 

After this deform continuously the blocks of size 2 so that they 
become diagonal (by Proposition~\ref{Pi} this is possible). The resulting 
quadruple is diagonal. It is a direct sum of $k+1$ quadruples 
$(\lambda _1,\lambda _2,\lambda _3,\lambda _4)$ and of $k$ quadruples 
$(\mu _1,\mu _2,\mu _3,\mu _4)$. It is unique up to conjugacy and can be 
reached by continuous deformation from any quadruple of ${\cal V}$. Hence, 
${\cal V}$ is connected.~$\Box$

\section{The variety ${\cal V}$ is not always connected\protect\label{secex}} 

We illustrate the title of the section by the following

\begin{ex} 
Consider the case $n=2$, $p=2$, the conjugacy classes $c_1$ and $c_2$ being 
diagonalizable, with 
eigenvalues $\pi ,2$ and $1-\pi,-1$, the conjugacy class $c_3$ consisting of 
the non-scalar matrices with eigenvalues $-1,-1$. Hence, $\kappa =2$ and the 
triple of conjugacy classes satisfies the conditions of 
Theorem~\ref{necessary} (to be checked directly). 
A priori ${\cal V}$ contains at least the following  
two components (denoted by ${\cal V}_1$ and ${\cal V}_2$). 
In ${\cal V}_1$ the triples 
of matrices $A_j$ equal (up to conjugacy)

\[ A_1=\left( \begin{array}{cc}\pi &1\\0&2\end{array}\right) ~~,~~
A_2=\left( \begin{array}{cc}1-\pi &0\\0&-1\end{array}\right) ~~,~~
A_3=\left( \begin{array}{cc}-1&-1\\0&-1\end{array}\right)\]
In ${\cal V}_2$ they equal (up to conjugacy)   

\[ A_1=\left( \begin{array}{cc}\pi &0\\1&2\end{array}\right) ~~,~~
A_2=\left( \begin{array}{cc}1-\pi &0\\0&-1\end{array}\right) ~~,~~
A_3=\left( \begin{array}{cc}-1&0\\-1&-1\end{array}\right)\]
The variety ${\cal V}$ contains no irreducible triples, see 
3) of Remarks~\ref{Gabber}. Hence, every triple from ${\cal V}$ 
is triangular up to conjugacy but 
not diagonal (otherwise $A_3$ must be scalar). Hence,  
${\cal V}={\cal V}_1\cup {\cal V}_2$. 

On the other hand,  
${\cal V}_1\cap {\cal V}_2=\emptyset$ because the eigenvalues with which 
each matrix acts on the invariant subspace are different for the two 
components. Hence ${\cal V}$ is disconnected.
\end{ex}

In the above example, however, the closure of ${\cal V}$ is connected. Indeed, 
consider the matrices 

\[ A_1'(\varepsilon )=\left( \begin{array}{cc}\pi &\varepsilon \\0&2
\end{array}\right) ~~,~~
A_2'(\varepsilon )=\left( \begin{array}{cc}1-\pi &0\\0&-1 
\end{array}\right) ~~,~~
A_3'(\varepsilon )=\left( \begin{array}{cc}-1&-\varepsilon \\0&-1\end{array}
\right)\]
where $\varepsilon \in ({\bf C},0)$. For $\varepsilon \neq 0$ this is a triple 
of matrices from ${\cal V}_1$, for $\varepsilon =0$ this is a triple from 
its closure (but not from ${\cal V}$ because $A_3$ is scalar). In the same 
way, for $\varepsilon \neq 0$ the matrices 

\[ A_1''(\varepsilon )=\left( \begin{array}{cc}\pi &0\\ \varepsilon &2
\end{array}\right) ~~,~~
A_2''(\varepsilon )=\left( \begin{array}{cc}1-\pi &0\\0&-1 
\end{array}\right) ~~,~~
A_3''(\varepsilon )=\left( \begin{array}{cc}-1&0\\-\varepsilon &-1\end{array}
\right)\]
belong to ${\cal V}_2$, for $\varepsilon =0$ they belong to its closure but 
not to ${\cal V}$ and one has $A_j'(0)=A_j''(0)$. 

In the above example the disconnectedness of ${\cal V}$ seems to result from 
the class $c_3$ not being closed. 
It would be interesting to prove or disprove that the closure of 
${\cal V}$ is always connected.

\section{Another example\protect\label{tex}}
\subsection{Description of the example}

{\em For values of $(-\kappa )$ arbitrarily big there exist examples when a 
component of ${\cal V}$ does not lie in the closure of the union of its 
components on which the centralizer is trivial, and is of dimension higher 
than the expected one.}

Indeed, consider the following example. 
Suppose that $p>3$ and that the $(p+1)$ conjugacy classes $c_j$ (or $C_j$) are 
diagonalizable, each MV being of the form $(m_j,1,1,\ldots ,1)$, 
$3\leq m_j\leq n-1$. Hence, 
$r_j=n-m_j$. Suppose that $r_1+\ldots +r_{p+1}=2n-2$. 

\begin{lm}
The $(p+1)$-tuple of conjugacy classes $c_j$ or $C_j$ like above satisfies the 
conditions of Theorem~\ref{necessary}.
\end{lm}

Indeed, one has $n_1=n-2$ and applying $\Psi$ once one obtains a $(p+1)$-tuple 
of conjugacy classes satisfying condition $(\omega _{n-2})$, see 
Remark~\ref{rd}.~~~~$\Box$

Denote by $\mu _j$ the 
eigenvalue of $A_j$ (or of $M_j$) of multiplicity $m_j$. Suppose that there 
holds the only non-genericity relation $\mu _1+\ldots +\mu _{p+1}=0~~(*)$ 
(resp. $\mu _1\ldots \mu _{p+1}=1$). 

\begin{rem}\label{R} 
There exists no irreducible $(p+1)$-tuple of matrices 
$A_j\in c_j$ whose sum is 0, see Remark~\ref{nongeneric}. Indeed,  
condition $(\delta _n)$ from Remark~\ref{nongeneric} does not hold -- 
set $b_j=\mu _j$ and recall that $(*)$ holds; 
then rk$(A_j-b_jI)=r_j$ and $r_1+\ldots +r_{p+1}=2n-2<2n$. A similar remark 
holds for matrices $M_j$ as well.
\end{rem}

Define the conjugacy classes $c_j^*\subset gl(n-2,{\bf C})$ and 
$c_j'\subset gl(n-1,{\bf C})$ (resp. 
$C_j^*\subset GL(n-2,{\bf C})$ and $C_j'\subset GL(n-1,{\bf C})$) as obtained 
from $c_j$ 
(resp. from $C_j$) by keeping the distinct eigenvalues the same and by 
decreasing the multiplicity of $\mu _j$ by 2 and by 1; the JNFs defined by the 
conjugacy classes $c_j^*$, $C_j^*$, $c_j'$ and $C_j'$ are diagonal. Hence, 
the sum (resp. the product) of all eigenvalues  of the classes 
$c_j^*$ and $c_j'$ (resp. $C_j^*$ and $C_j'$) counted with the 
multiplicities equals 0 (resp. 1). 

Condition $(\omega _{n-2})$ holds for the classes $c_j^*$ or $C_j^*$ while 
condition $(\omega _{n-1})$ holds for the classes $c_j'$ or $C_j'$. By 
Theorem~2 from \cite{Ko4} (we need $p>3$ to apply it), the DSP is 
solvable for the classes $c_j^*$ (resp. $C_j^*$) and $c_j'$ (resp. $C_j'$). 
Denote by $H_j\in c_j^*$ (resp. $H_j\in C_j^*$) and $G_j\in c_j'$ 
(resp. $G_j\in C_j'$) matrices with sum equal to 0 (resp. with product equal 
to $I$) whose $(p+1)$-tuple is irreducible. 

\begin{prop}\label{H}
1) There exist $(p+1)$-tuples of 
$n\times n$-matrices with trivial centralizers, whose sums equal 0 (or whose 
products equal $I$) and blocked as follows: 
\[ A_j~{\rm (or~}M_j{\rm )}=\left( \begin{array}{ccc}H_j&R_j&Q_j\\0&\mu _j&0\\
0&0&\mu _j\end{array}\right) ~~{\rm or}~~A_j~{\rm (or~}M_j{\rm)}=
\left( \begin{array}{ccc}\mu _j&0&T_j\\
0&\mu _j&S_j\\
0&0&H_j\end{array}\right) ~.\] 

2) Any $(p+1)$-tuple with trivial centralizer of matrices $A_j\in c_j$ or 
$M_j\in C_j$ is up to conjugacy block upper-triangular, with all diagonal 
blocks but one being equal, of size one, the restriction of $A_j$ or $M_j$ 
to such a 
block being equal to $\mu _j$. The different block is first or last on the 
diagonal. The number of diagonal blocks is $\geq 3$.
\end{prop}

The proposition is proved in the next subsection. 

Consider the 
stratum ${\cal U}\subset {\cal V}$ of $(p+1)$-tuples of 
matrices which up to conjugacy are of the form 
$\tilde{G}_j=\left( \begin{array}{cc}G_j&0\\0&\mu _j\end{array}\right)$. 

\begin{lm}
A point of the stratum ${\cal U}$ does not belong to the closure of any of 
the strata on which the centralizer is trivial.
\end{lm}

Indeed, the matrix algebra 
generated by the $(p+1)$-tuples of 
matrices defined by a point of the stratum ${\cal U}$ 
contains a matrix with distinct eigenvalues (the 
$(p+1)$-tuple of matrices $G_j$ is 
irreducible and defines a representation not equivalent 
to $(\mu_1,\ldots ,\mu _{p+1})$) while each 
matrix from an algebra defined by a point of 
any stratum of ${\cal V}$ where the centralizer is trivial has an 
eigenvalue of multiplicity $\geq 2$, see 2) of Proposition~\ref{H}.~~~~~$\Box$

\begin{prop}\label{UW}
1) One has dim${\cal U}=3(n-1)^2+1-\sum _{j=1}^{p+1}r_j^2$.

2) The dimension of each of the two strata of ${\cal V}$ (denoted by 
${\cal W}_1$, ${\cal W}_2$) of the $(p+1)$-tuples of matrices 
which up to conjugacy are like the ones from 1) of Proposition~\ref{H} 
equals $3(n-1)^2-\sum _{j=1}^{p+1}r_j^2=$dim${\cal U}-1$.
\end{prop}

The proposition is proved in Subsection~\ref{prUW}. It implies that 
dim${\cal U}>$dim${\cal W}_i$, i.e. dim${\cal U}$ is greater than the expected 
dimension.

\subsection{Proof of Proposition~\protect\ref{H}}

We prove the proposition only in the 
case of matrices $A_j$ and for the left $(p+1)$-tuple of matrices given in 
1) of the proposition leaving for the 
reader the proof in the other cases -- it can be performed in a 
similar way. We part prove 1) in $1^0$ -- $2^0$ and part 2) in 
$3^0$ -- $4^0$. 

$1^0$. Denote by $H$ an irreducible $(p+1)$-tuple of matrices like in the 
proposition as well as the representation defined by it and by $\mu$ the 
$(p+1)$-tuple $(\mu _1,\ldots ,\mu _{p+1})$.  

One has 
$\delta :=$dim~Ext$^1(H,\mu )=$dim~Ext$^1(\mu ,H)=2$. Indeed,  
$\delta =$dim$({\cal L}/{\cal N})$ where 

\[ {\cal L}=\{ (L_1,\ldots ,L_{p+1})~|~L_j=(H_j-\mu _jI)X_j~,~
L_1+\ldots +L_{p+1}=0\} ~,\]
\[ {\cal N}=\{ ((H_1-\mu _1I)X,\ldots ,(H_{p+1}-\mu _{p+1}I)X)\} \]
where the matrices $X_j$ and $X$ are $(n-2)\times 1$. 

The dimension of the space of matrices of the form $(H_j-\mu _jI)Y=$ 
(where $Y$ is $(n-2)\times 1$) equals $r_j$. The condition 
$L_1+\ldots +L_{p+1}=0$ is equivalent to $n-2$ linearly independent conditions 
(their linear independence follows easily from the fact that the 
representation $H$ is irreducible and not equivalent to the one-dimensional 
representation $\mu$). Hence, dim${\cal L}=r_1+\ldots +r_{p+1}-n+2=n$. 
The same kind of argument shows that dim${\cal N}=n-2$ which 
implies that $\delta =2$.  

$2^0$. It follows from $1^0$ that one can construct two linearly indepent 
$(p+1)$-tuples of $(n-2)\times 1$-matrices belonging to the space 
${\cal L}/{\cal N}$ -- these are the $(p+1)$-tuples of matrices $R_j$ and 
$Q_j$. Show that the centralizer of the thus constructed $(p+1)$-tuple of 
matrices is trivial. Denote a matrix from the centralizer by 
$Z=\left( \begin{array}{ccc}K&B&C\\D&e&f\\U&v&w\end{array}\right)$ where 
$K$ is of the size of $H_j$ etc. 

The commutation relations imply 
$UH_j-\mu _jU=0$, $j=1,\ldots ,p+1$. It follows from $H$ and $\mu$ being 
non-equivalent and 
$H$ being irreducible that $U=0$. But then $DH_j-\mu _jD=0$, 
$j=1,\ldots ,p+1$, and in the same way one obtains $D=0$. 

Hence, one has $[K,H_j]=0$, $j=1,\ldots ,p+1$, which implies that 
$K=\alpha I$, $\alpha \in {\bf C}$ (recall that $H$ is irreducible). 

This in turn implies that 

\[ (H_j-\mu _jI)B+(e-\alpha )R_j+vQ_j=0~~,~~
(H_j-\mu _jI)C+fR_j+(w-\alpha )Q_j=0~.\]
The definition of the matrices $R_j$ and $Q_j$ implies that $B=C=0$, 
$v=f=0$, $e=w=\alpha$. Hence, the centralizer is trivial. This proves 1) of 
the proposition.

$3^0$. Prove 2). Recall that there exists no irreducible $(p+1)$-tuple of 
matrices $A_j\in c_j$ whose sum is 0 (Remark~\ref{R}) and that $(*)$ 
is the only non-genericity relation satisfied by the eigenvalues.
Hence, every $(p+1)$-tuple of matrices $A_j\in c_j$ whose sum is 0 is up to 
conjugacy block upper-triangular and all diagonal blocks but one (denoted 
by $D$) are of size 1 
and the restrictions of $A_j$ to them equal $\mu _j$. (The block $D$  
can also be of size 1 but in this case $A_j|_D\neq \mu _j$.)

$4^0$. If the first and the last diagonal blocks are equal, 
then the centralizer of the $(p+1)$-tuple is non-trivial -- it contains the 
matrix $E_{1,n}$. So assume that the first diagonal block is different from 
all others (the case when this is the last block can be treated in a similar 
way). If there are only two diagonal blocks ($D$ of size $n-1$ and $\mu$ of 
size 1), then one has dim~Ext$^1(D,\mu )=0$ (this follows from 
$r_1+\ldots +r_{p+1}=2n-2$) and, hence, the representation defined by the 
matrices $A_j$ is a direct sum. Hence, if the centralizer of the 
$(p+1)$-tuple is trivial, then there are at least three diagonal 
blocks.~~~~~$\Box$

\subsection{Proof of Proposition \protect\ref{UW}\protect\label{prUW}}

We prove the proposition only for matrices $A_j$, for matrices $M_j$ the 
proof is similar.

$1^0$. Prove 1). 
The dimension of the variety of irreducible $(p+1)$-tuples of matrices 
$G_j\in c_j'$ with zero sum equals 
$u'=(\sum _{j=1}^{p+1}d(c_j'))-((n-1)^2-1)$, see Remarks~\ref{expected}. 
Hence, the dimension of block-diagonal matrices $\tilde{G}_j$ whose 
sum is 0 equals $u'$. To obtain dim${\cal U}$ one has to add to $u'$ the 
dimension of a 
transversal $T$ at $I$ to the subgroup of $GL(n,{\bf C})$ of 
block-diagonal matrices with diagonal blocks of sizes $n-1$ and 1, the only 
ones conjugation with which preserves the block-diagonal form of the 
$(p+1)$-tuple. One has $d(c_j')=(n-1)^2-r_j-(n-1-r_j)^2=(2n-3)r_j-r_j^2$ (see 
Remark~\ref{rd}) and dim$T=2n-2$. Hence, 

\[ {\rm dim}{\cal U}=(2n-3)\sum _{j=1}^{p+1}r_j-\sum _{j=1}^{p+1}r_j^2-
((n-1)^2-1)+2n-2=3(n-1)^2+1-\sum _{j=1}^{p+1}r_j^2~.\] 
 
$2^0$. Prove 2). The dimension of the 
variety of $(p+1)$-tuples of matrices $H_j\in c_j^*$ 
whose sum is 0 equals $u^*=(\sum _{j=1}^{p+1}d(c_j^*))-((n-2)^2-1)$ 
(computed like $u'$, by changing $n-1$ to $n-2$). Hence, this is 
the dimension of the variety of $n\times n$-matrices which are block-diagonal, 
with diagonal blocks equal to $H_j$, $\mu _j$, $\mu _j$. Note that 
$d(c_j^*)=(n-2)^2-r_j-(n-2-r_j)^2=(2n-5)r_j-r_j^2$ (Remark~\ref{rd}).

$3^0$. The dimension of each of the two varieties of $(p+1)$-tuples of 
matrices like in 1) of Proposition~\ref{H} equals 
$u^*+2\sum _{j=1}^{p+1}r_j-2(n-2)=u^*+2n$. Indeed, each of 
the matrices $Q_j$, $R_j$ or $T_j$, $S_j$ belongs to a linear space of 
dimension $r_j$ (this is the image of the linear operator 
$(.)\mapsto (H_j-\mu _jI)(.)$ or $(.)\mapsto (.)(H_j-\mu _jI)$ acting on 
${\bf C}^{n-2}$). One has to subtract $2(n-2)$ because $\sum R_j=\sum Q_j=0$ 
and $\sum T_j=\sum S_j=0$.

$4^0$. One can consider the matrices from 1) of Proposition~\ref{H} like block 
upper-triangular, with two diagonal blocks the lower of which is of size 2 
and is scalar. The subgroup of $GL(n,{\bf C})$ conjugation with which 
preserves this form is the subgroup of block upper-triangular matrices 
with diagonal blocks of sizes $n-2$ and 2. Hence, a transversal at $I$ to it 
is of dimension $2(n-2)$ and one has 

\[ {\rm dim}{\cal W}_i=u^*+2n+2(n-2)=u^*+4n-4=\]
\[ =(2n-5)\sum _{j=1}^{p+1}r_j-\sum _{j=1}^{p+1}r_j^2-((n-2)^2-1)+4n-4=
3(n-1)^2-\sum _{j=1}^{p+1}r_j^2~.~~~~~\Box\]

Author's address: Universit\'e de Nice, Laboratoire de Math\'ematiques, 
Parc Valrose, 06108 Nice Cedex 2, France

e-mail: kostov@math.unice.fr
\end{document}